\newtheorem{defn}{Definition}[section]
\newtheorem{lem}[defn]{Lemma}
\newtheorem{thm}[defn]{Theorem}
\newtheorem{prop}[defn]{Proposition}

\newtheorem{rem}[defn]{Remark}

\newtheorem{assu}[defn]{Assumption}

\newtheorem{notdef}[defn]{Notation - Definition}

\newcommand{\A}{{\bf A}}
\newcommand{\eps}{\epsilon}
\newcommand{\pp}{{\bf P}}
\newcommand{\pgl}{{\bf PGL}}
\newcommand{\pabc}{\pp_{a,b,c}}
\newcommand{\tpabc}{\tilde{\pp}_{a,b,c}}
\newcommand{\op}{{\cal O}_{\pp^1}}
\newcommand{\proof}{{\bf Proof:}\,\,}
\newcommand{\inv}{^{-1}}
\newcommand{\Pic}{{\rm Pic}}
\newcommand{\qed}{$\quad \diamond$\par\smallskip}
\newcommand{\OO}{{\cal O}}
\newcommand{\OX}{{\cal O}_X}
\newcommand{\OY}{{\cal O}_Y}

\newcommand{\OP}{{\cal O}_{\pp^1}}

\newcommand{\OSi}{{\cal O}_{\Si}}
\newcommand{\HH}[1]{{\rm H}^{#1}}
\newcommand{\tF}{\tilde{F}}
\newcommand{\Si}{\Sigma}

\newcommand{\Proj}{{\rm Proj}}
\newcommand{\Z}{{\bf Z}}

\title{Irregular canonical double surfaces}
\author{Margarida  Mendes Lopes \quad --- \quad Rita
Pardini
\thanks{The first author is a member of CAUL and the second author is a
member of GNSAGA of CNR. The present collaboration takes place within
the framework of contract CHRXCT940557 of the HCM program of EU.}}
\date{}
\documentstyle{article}
\begin{document}
\maketitle
\setcounter{section}{-1}
\def\theequation{\thesection.\arabic{equation}}

\section{Introduction}
\setcounter{defn}{0}
\setcounter{equation}{0}

Let $X$ be a minimal surface of general type of geometric genus $p_g$,
let $\Si\subset\pp^{p_g-1}$ be the canonical image of $X$ and let
$\phi:X\to\Si$ be the canonical map.
If $\Si$ is a surface but $\phi$ is not birational, then by theorem $3.1$
of
\cite{bea} either i)
$p_g(\Si)=0$ or ii) $\Si$ is the canonical image of a surface of general
type $S$ whose canonical map is birational (and then, of course,
$p_g(\Si)=p_g$).
 Recall (cf. \cite{bea}, thm.
$5.5$ or \cite{ha}, page $44$) that the Castelnuovo inequality $K^2\ge
3p_g-7$ holds for surfaces of general type with birational canonical map;
so in case ii)
  $S$ satisfies $K^2_S\ge 3p_g-7$, and $K^2_X\ge 6p_g-14$, with equality
holding if and only if the canonical system of $X$ is base point free
and the minimal resolution $S$ of $\Si$ is on the
Castelnuovo line $K^2_S=3p_g-7$ (cf. prop \ref{harris} and proof).

Case ii) of the theorem quoted above was thought to be impossible for a
long time. In fact only very few examples are known, and all but
one series due to Beauville (cf. section \ref{esempi}, example $4$)
have bounded invariants. The examples in this infinite series satisfy:
$K^2_X=6p_g-14$, $\deg
\phi=2$ and $q(X):=h^0(X,\Omega^1_X)=2$.

The main purpose of this paper is to show that these are almost the
only examples satisfying $K^2_X=6p_g-14$ and $q(X)\geq 2$. Therefore we
make the following:
\begin{assu}\label{ipointro}
Let $X$ be a minimal surface of general type  of geometric genus $p_g$,
with $K^2_X=6p_g-14$, $q(X)>0$ and $K_X$ ample, let
$\Si\subset\pp^{p_g-1}$ be the canonical image of $X$ and let
$\phi:X\to\Si$ be the canonical map: assume that $\Si$ is a canonical
surface and that
$\phi$ is not birational. Moreover, if $p_g=4,5,7$  assume that $\Si$ is
isomorphic to a divisor with at most rational double points in a
$\pp^2$-bundle over $\pp^1$, such that the fibres $F$ of the projection
$\Si\to\pp^1$ are plane quartics.
\end{assu}

If the above assumption is satisfied, then the minimal
desingularization $S$ of $\Si$ satisfies $K^2_S=3p_g-7$.
Surfaces with these numerical invariants have been described by Ashikaga
and Konno in \cite{ak}: for $p_g=6$ or $p_g\ge 8$, those with birational
canonical map are isomorphic to a divisor with at most rational double
points in a $\pp^2$-bundle over $\pp^1$, such that the fibres $F$ of the
projection $\Si\to\pp^1$ are plane quartics. (This accounts for the
somewhat funny-looking final part of assumption \ref{ipointro}.)
We divide
the surfaces $X$ in types $I$ and $I\!I$, according to whether, for a
generic fibre $F$,
$\phi^*F$ is connected or not.
Surfaces of type $I$ are the ``general case'', and, if $q(X)\ge 2$, they
correspond to Beauville's examples:
\begin{thm}\label{mainone}
Assume that \ref{ipointro} holds, that $q(X)\ge 2$ and that $X$ is of type
$I$: then  $p_g(X)\equiv 1\,\, (\mbox{mod}\,\,\, 3)$, $q(X)=2$, the
Albanese surface $A$ of $X$ has an irreducible principal polarization,
and $X$ can be constructed as in example
$4$ of section \ref{esempi}, with $n=(p_g(X)+3)/4$.
\end{thm}
Let us remark that we do not know any example of type $I$ surfaces with
$q(X)=1$. To establish whether such surfaces exist, and in case of
existence whether they have bounded invariants is an interesting
problem.

 On the other hand, surfaces of type
$I\!I$ should be regarded as exceptions, and can be described completely:
\begin{thm}\label{maintwo}
Assume that \ref{ipointro} holds and that $X$ is of type
$I\!I$, let $X\to B\to\pp^1$ be the Stein factorization of the pencil
$\phi^*|F|$ and let $g$ be the genus of $B$: then there exist integers
$0\le a\le b \le c$ with $c\le g$ and $a+b+c=p_g-3$ such that $\Si$ is
isomorphic to a divisor in $\pabc:=\Proj(\OP(a)\oplus\OP(b)\oplus\OP(c))$
with the following properties: i)
$\Si$
 is linearly equivalent to
$4T-(a+b+c-2)L$, where $T$ is the tautological hyperplane section and $L$
is the fibre of $\pabc$ (and $F=L|_{\Si}$), ii) the pencil $|F|$ on $\Si$
has precisely
$2g+2$ double fibres, iii) the only singularities of
$\Si$ are nodes and $\Si$ is smooth outside the double fibres of $|F|$.
 The
double fibres of $|F|$ occur at the branch points of $B$ and each contains
$8$ nodes .

Conversely, given integers $0\le a\le b\le c$ and $g$, with $c\le g$, if
$\Si\subset \pabc$ is a divisor satisfying conditions i),ii),iii)
above, then
$\Si$ has
$16g+16$ nodes and there exists a double cover $\phi:X\to\Si$ branched
over the nodes such that $X$ is a surface of type $I\!I$ and $\phi$ is the
canonical map of $X$.

The numerical possibilities for the invariants of $X$ are the following:
\smallskip

\noindent a) $p_g=3g+3$,\quad $q=g$,\quad\quad $a=b=c=g$,\quad $0<g\le 26$;

\noindent b) $p_g=3g+2$,\quad $q=g+1$, \quad $a=b=g-1$, $c=g$, \quad
$0<g\le 17$;

\noindent c) $p_g=3g+1$,\quad $q=g+2$,\,\,\,\quad $a=b=g-1$, $c=g$ or
 $a=g-2$,
$b=c=g$, \quad
$0<g\le 8$.
\end{thm}

In \cite{bea} it is also proven that if the canonical map $\phi:X\to\Si$
is not birational and $\Si$ is a canonically embedded surface then
$\deg\phi\le 3$ for
$\chi(\OX)\ge 14$, and if $\deg\phi=3$ then
$q(X)\le 3$. Thm.\ref{mainone} and \ref{maintwo} imply in particular
that the irregularity
$q(X)$ is also bounded under assumption \ref{ipointro}. It would be
interesting to know whether $q(X)$ is bounded in general for $\deg\phi=2$.
Another interesting problem is to study {\em regular} surfaces $X$ such
that the canonical map is not birational and the canonical image is a
canonically embedded surface: only very few examples are known (cf.
section \ref{esempi}) and, lacking the information given by the Albanese
map, their structure is quite mysterious even when the invariants satisfy
the ``minimal'' relation $K^2_X=6p_g-14$.

The paper is organized as follows: in section $1$ we set the notation
and recall some facts on double cover that will be used later. In section
$2$ we describe the general set-up and establish various facts about $X$,
$S$ and
$\Si$. In particular we study the structure of degenerate fibres of
$\phi^*F$ both for type $I$ and type $I\!I$. In section $3$ we describe the
construction of all the examples known to us of surfaces of general type
with $2$-$1$ canonical map onto a canonical surface. In section
$4$ we look at the surfaces of type $I$ with $q(X)\ge 2$ and we show, using
a fine analysis involving the Albanese map and the Prym variety of
$\phi^*F\to F$ for general $F$, that these are exactly Beauville's
examples. In section
$5$ we describe the surfaces of type $I\!I$ in detail and determine the
possible ranges for their invariants. Section $6$ contains a computation
with Macaulay that shows that  example $3$ of section $3$ actually exists.

{\em Acknowledgements:} We are indebted to several people for useful
conversations; we would like to mention in particular C. Birkenhake, F.
Catanese, M. Manetti, G. Ottaviani, N. Shepherd-Barron, A. Verra.

\section {Notation and conventions}\label{notation}
\setcounter{defn}{0}
\setcounter{equation}{0}

All varieties are normal projective varieties over the complex numbers.
The $n$-dimensional projective space is denoted by $\pp^n$, and its
group of automorphisms by $\pgl(n)$. As usual,
$\OY$ is the structure sheaf of the variety
$Y$,
$\HH{i}(Y,{\cal F})$ is the $i$-th cohomology group of a sheaf ${\cal
F}$ on $Y$, and $h^i(Y,{\cal F})$ is the dimension of $\HH{i}(Y,{\cal
F})$; for a line bundle $M$ on $Y$, we denote by $|M|$ the complete
linear system $\pp(\HH{0}(Y,M))$. When dealing with smooth varieties, we
do not distinguish between line bundles and divisors. If
$Y$ is smooth, then
$K_Y$ denotes a canonical divisor and $\Pic(Y)$ the Picard group of
$Y$. If $Y$ is a surface, then
$p_g(Y)=h^0(Y,K_Y)$ is the {\em geometric genus} and
$q(Y)=h^1(Y,\OO_Y)$ is the {\em irregularity}, $K_Y^2$ is the
self-intersection of the canonical divisor; we denote by
$\chi(Y)=1-q(Y)+p_g(Y)$  the Euler characteristic of $\OY$ and by $c_2(Y)$
the second Chern class of the tangent bundle of $Y$, or, which is the
same, the topological Euler characteristic of $Y$. A surface
$Y$ is said to be {\em irregular} if $q(Y)\ne 0$. The intersection
number of two divisors
$C$,
$D$ on a smooth surface is denoted simply by $CD$, linear equivalence is
denoted by $\equiv$.  A {\em node} of a surface is a double point of type
$A_1$, namely a hypersurface singularity that in suitable local analytic
coordinates is defined by the equation
$x^2+y^2+z^2=0$.

A {\em double cover} is a finite map $f:X\to Y$ of degree
$2$ between normal projective varieties; we denote by $i:X\to X$ the
involution that interchanges the two points of a generic fibre of
$f$.
In this paper we
will need to consider only the following two cases: a) both $X$ and $Y$ are
smooth, and b) $X$ is a smooth surface, $Y$ is normal and $f$ is
unramified in codimension $1$.

In case a), $f$ is a flat map and $f_*\OX$
splits under the action of $i$ as $O_Y\oplus {\cal L}\inv$, where ${\cal
L}$ is a line bundle and $i$ acts on
${\cal L}\inv$ as multiplication by $-1$. The branch locus of $f$ is a
smooth divisor
$B\equiv 2{\cal L}$, the ramification locus is a divisor $R\equiv
f^*{\cal L}$ and one has:
\begin{eqnarray}\label{formuledoppi}
K_X=f^*(K_Y+{\cal L})\quad K_X^2=2(K_Y+{\cal L})^2\quad f_*K_X=K_Y\oplus
K_Y+{\cal L}\\ h^i(\OX)=h^i(\OY)+h^i({\cal L}),\quad i=1,\ldots \dim Y
\end{eqnarray}
(Actually, the above formulas also hold if $Y$ is a surface with rational
dpoble points and $B$ is a smooth divisor containing no singulairities
of $Y$).
 The cover
$\phi:X\to Y$ can be reconstructed from $Y$, ${\cal L}$, $B$ as follows.
Let
$p:{\cal L}\to Y$ be the projection, let $w$ be the tautological section of
$p^*{\cal L}$ and let
$\sigma\in
\HH{0}(Y,{\cal L}^2)$ be a section vanishing on $B$: the zero locus in
${\cal L}$ of the section $w^2-p^*\sigma$ of $p^*{\cal L}^2$, together
with the restriction of the map $p$, is a double cover of $Y$ isomorphic to
$\phi:X\to Y$. Moreover, it is clear that, given a line bundle ${\cal L}$
on
$Y$ and a divisor $B$ in the linear system $|{\cal L}^2|$, the above
construction yields a finite degree $2$ map $\phi:X\to Y$.

A {\em linearization} of a line bundle $N$  on $X$ is an
involution $i_N:N\to N$  that lifts the
involution
$i:X\to X$.  If $N$ is a linearized line bundle, we say that $\sigma\in
\HH{0}(X, N)$ is {\em even\,} if ${i_N}_*\sigma=\sigma$ and {\em odd\,} if
${i_N}_*\sigma=-\sigma$. A divisor defined by an even (odd) section is
called {\em symmetric\,}({\em antisymmetric}).
The canonical bundle
$K_X$ and the pull-backs of line bundles from $Y$ have natural
linearizations: in these cases, unless otherwise stated, we  consider the
natural linearizations.

Consider now case b): the singularities of $Y$ are nodes, that  are the
images of the fixed points of $i$. If $\nu$ is the number of nodes of
$Y$, then one has (see \cite{becky} (0.6)):
\begin{equation}\label{nodi}
\chi(\OX)=2\chi(\OY)-\frac{1}{4}\nu
\end{equation}
 A set
$J$ of nodes on a normal surface $Y$ is said to be {\em even} if there
exists a double cover $\phi:X\to Y$ branched precisely over $J$.
\begin{prop}\label{criterio}
Let $W$ be a smooth $3$-fold, let $Y\subset W$ be a divisor whose only
singularities are nodes;  if there exists a divisors $D$, $D'$  in
$W$ such that  $D\equiv 2D'$ and $D$ restricted to $Y$ is
equal to
$2C$, where $C$ is a curve passing though all the nodes of $Y$, then the
nodes of $Y$ are an even set.
\end{prop}
\proof
Denote by $\eta:\hat{W}\to W$ the blow-up at the nodes of $Y$, by
$\epsilon:\hat{Y}\to Y$ its restriction to the strict transform
$\hat{Y}$ of
$Y$, by $\hat{E_i}$, $E_i$ the exceptional divisors of
$\eta$ and $\epsilon$ respectively, and by
$\hat{D}$,
$\hat{C}$ the strict transforms of
$D$ on
$\hat{W}$ and of $C$ on $\hat{W}$. One has the following linear
equivalence on $\hat{W}$:
$2\eta^*D'\equiv\eta^*D=\hat{D}+\sum E_i$, which
restricts to
$2\epsilon^*D'\equiv\epsilon^*D=2\hat{C}+\sum E_i$. So
$\sum E_i\equiv 2(\epsilon^*D'-\hat{C})$, and there exists a
smooth double cover $g:\hat{X}\to \hat{Y}$ branched over $\sum E_i$; the
ramification divisor of $g$ is a union of disjoint $-1$ curves that can
be contracted to yield $f:X\to Y$ branched over the nodes of $Y$.
\qed

\section{The set-up}\label{setup}
\setcounter{defn}{0}
\setcounter{equation}{0}
The notations and the assumptions introduced in this section will be
maintained throughout all the paper. We start by making the following:
\begin{assu}\label{ipotesi}
Let $X$ be a minimal surface of general type  of geometric genus $p_g$,
with $K^2_X=6p_g-14$, $q(X)>0$ and $K_X$ ample, let
$\Si\subset\pp^{p_g-1}$ be the canonical image of $X$ and let
$\phi:X\to\Si$ be the canonical map: assume that $\Si$ is a
canonical surface and that
$\phi$ is not birational. Moreover, if $p_g=4,5,7$  assume that $\Si$ is
isomorphic to a divisor with at most rational double points in a
$\pp^2$-bundle over $\pp^1$, such that the fibres $F$ of the projection
$\Si\to\pp^1$ are plane quartics.
\end{assu}
Let now $A$ be the Albanese variety of $X$, let $x_0\in X$ be a fixed
point of
$i$, let
$\alpha:X\to A$ be the Albanese map with base point $x_0$ and let
$K=A/<-1>$ be the Kummer variety of $A$. Since $\Sigma$ is regular,
the involution $i$ on $X$ induces on $A$ the multiplication by $-1$,
and there is an induced map $f:\Si\to K$. Thus we have the following
{\em basic commutative diagram}, where $q:A\to K$ is the natural
projection:
\begin{equation}\label{diagram}
\begin{array}{rcccl}
\phantom{1} & X &\stackrel{\alpha}{\rightarrow} & A & \phantom{1} \\
\scriptstyle{\phi}\!\!\!\!\!\! & \downarrow & \phantom{1} & \downarrow
& \!\!\!\!\!\!
\scriptstyle{q}
\\
\phantom{1} & \Si & \stackrel{f}{\rightarrow} & K & \phantom{1}
\end{array}
\end{equation}
Remark that, since $\phi$ is finite, $\phi:X\to\Si$ is obtained from
$q$ by base change and normalization.

Assuming $K^2_X=6p_g-14$ is equivalent to considering the lowest possible
value of $K^2$ in the above situation, as it appears from the next
proposition and its proof. In order to state it we introduce the
following
\begin{notdef}\label{castelnuovo}
Let $0\le a\le b \le c$ be integers: we write
$\pabc =
\mbox{Proj}(\op(a)\oplus\op(b)\oplus\op(c))$, and denote by $T$ the
tautological hyperplane section and by $L$ the fibre of
$\pabc$. We define a {\rm Castelnuovo surface of type $(a,b,c)$} to be a
divisor $\Si$ in $\pabc$ linearly equivalent to $4T-(a+b+c-2)L$ with
at most rational double points as singularities. Notice that $T$
restricts to the canonical divisor of $\Si$ and that the minimal
desingularization $S$ of $\Si$ satisfies: $K_S^2=3(a+b+c)+2=3p_g(S)-7$,
$q(S)=0$.
\end{notdef}
\begin{prop}\label{harris}
Assume that \ref{ipotesi} holds: then $K_X$ is base point free, the
degree of $\phi$ is equal to $2$,  $K^2_S=3p_g-7$, the only
singularities of $\Si$ are nodes,
and there exist integers
$0\le a\le b
\le c$ with
$a+b+c=p_g-3$ such that
$\Si$ is a Castelnuovo surface of type $(a,b,c)$.
\end{prop}
\proof
Write $K_X=M+Z$, where $M$ is the moving part and $Z$ is the fixed part
of $K_X$ and denote by $d$ the degree of
$\Si\subset\pp^{p_g-1}$. By \cite{ha} page $44$, one has $d\ge 3p_g-7$
and thus:
\begin{eqnarray*}\label{ksquare}
6p_g-14=K_X^2=K_XM+K_XZ\ge K_XM= M^2+MZ\ge M^2\ge \\
(\deg\phi)d\ge
(\deg\phi)(3p_g-7).
\end{eqnarray*}
The first and the second inequality are consequences of
the fact that $K_X$ and $M$ are nef. It follows that $\deg\phi$ is equal
to
$2$ and all the above inequalities are equalities. In particular, one has
$K_XZ=0$ and $M^2=2d$, and so $Z$ is empty (recall that
$K_X$ is ample) and $K_X=M$ is base point free. The surface $\Si$
satisfies $d=3p_g-7$: so it is Castelnuovo variety in the sense of
\cite{ha}, and in particular (cf. \cite{ha}, page $66$) it is
projectively normal, and therefore normal. So $\phi:X\to \Si$ is a double
cover as defined in section \ref{notation}, the only singularities of
$\Si$ are nodes, and therefore $\Si$ is the canonical model of $S$.
Castelnuovo varieties are classified in
\cite{ha} (apart from a small mistake corrected in \cite{mi}), and the
complete list of those that are  canonical surfaces is given in thm.
$1.5$ of \cite{ak}, where they are studied in detail. In particular, for
$p_g=6$ or
$p_g\ge 8$ $\Si$ is a Castelnuovo surface  of type $(a,b,c)$  for some
$0\le a\le b\le c$ with $a+b+c+3=p_g$.
\qed

Denote by $|F|$ the pencil on $\Si$ induced by the projection
$\pabc\to\pp^1$ and by $|\tF|$ the pull-back $\phi^*|F|$.
Surfaces $X$ as in assumption \ref{ipotesi} fall into two types
according to the nature of $|\tilde{F}|$:
\begin{defn}
We say that a surface $X$ as in assumption \ref{ipotesi} is of type
$I$ if $|\tilde{F}|$ is irreducible, and of type $I\!I$ if $|\tilde{F}|$ is
reducible.
\end{defn}
Remark that if $X$ is of type $I$ then $|\tF|$ is a linear pencil of genus
$5$. We will show later (proposition \ref{invariantsII}) that, if $X$ is
of type
$I\!I$, then there are only a finite
number of numerical possibilities for the invariants of $X$, so that  one
should think of surfaces of type $I$ as of the ``general case''.
 We close this section by stating some general facts about $X$ and
$\Sigma$.
\begin{prop}\label{basic} Assume that assumption \ref{ipotesi} holds. Then
the map $\phi:X\to\Sigma$ is ramified precisely over the singular
locus of $\Sigma$, which consists of $4(1+p_g(X)+q(X))$ nodes.
\end{prop}
\proof
We have already remarked in prop. \ref{harris} and its proof that
$\phi:X\to \Si$ is a double cover.  Moreover,
$\phi$ is necessarily unramified in codimension $1$, since otherwise $K_X$
would have a fixed part. So the number of nodes of $\Si$ can be computed by
means of formula
\ref{nodi}. \qed

\begin{prop}\label{nodes}
Every irreducible component of a fibre of the pencil $|F|$ on $\Si$ has
multiplicity $\le 2$, and every double fibre contains $8$ singular
points of $\Si$. Moreover, if
$X$ is of type
$I$, then
$|F|$ and $\phi^*|F|$ have no multiple fibres.
\end{prop}
\proof  The fibres  $F$  on $\Si$ are (possibly singular) plane
quartics.
Since the only singularities of $\Si$ are nodes and the fibres $F$
are the restriction to $\Si$ of smooth Cartier divisors on the smooth
threefold $\pabc$, the fibres of
$|F|$ have double points at the singular points of
$\Si$. So, if $C$ is a component of a fibre of multiplicity
$m>2$, then $C$ does not contain any singular point of $\Si$.
$C$ is necessarily a line and thus,  if
$C'$ denotes the strict transform of
$C$ on $S$, one has: $K_S C'=1$, $C'^2=-3$, $C'(\epsilon^*F-mC')=4-m$.
(The last equality is a consequence of the fact that $C'$ contains no
singular point of $\Si$ and $F$ is a plane quartic.) So
it follows:
$0=C'\epsilon^*F=mC'^2+C'(\epsilon^*F-mC')=-3m+4-m=4(1-m)$, a
contradiction  since $m>1$. Let now $2C$ be a double fibre on $\Si$,
with $C$  an irreducible plane conic. Let $P_1,\ldots P_k$ be the nodes of
$\Si$ that lie on $C$, let $E_1,\ldots E_k$ be the corresponding
$-2$-curves on $S$ and let $C'$ be the strict transform of $C$ on $S$.
The pull-back of the  fibre $2C$ to $S$ is
$\epsilon^*(2C)=2C'+E_1+\ldots+E_k$, and
$C'E_i=1$, for $i=1,\ldots k$. If
$C$ is irreducible, then we have: $K_SC'=2$, $C'^2=-4$,
$0=C'\epsilon^*F=C'(2C'+E_1+\ldots+E_k)=-8+k$, $k=8$. If $C$ consists
of a pair of distinct lines, then a similar computation shows that each
line contains $4$ nodes of $\Si$.

Assume now that the surface $X$ is of type $I$, and suppose that $2C$
is a double fibre on $\Si$. The pull-back of $2C$ to $X$ is a double
fibre
$\tF=2D$, where $D$ is either a smooth hyperelliptic curve of genus
$3$, or the sum of two smooth elliptic curves meeting transversely at
$2$ points, according to whether $C$ is irreducible or not. In both
cases, $D$ is a curve of arithmetic genus $3$. Tensoring with
$K_X+\tF$ the decomposition sequence  $0\to\OO_D(-D)\to
\OO_{\tF}\to\OO_D\to 0$
and taking global sections, one obtains the following
exact sequence:
$$0\to \HH{0}(D,K_D)\to \HH{0}(\tF,K_{\tF})\to
\HH{0}(D,\OO_D(K_X+\tF)).$$
By the Ramanujan vanishing theorem, one has $\HH{1}(X,-\tF)=0$ and, as a
consequence,
$h^0(\tF,\OO_{\tF})=1$. As the arithmetic genus of $\tF$ is equal to $5$,
this implies that $h^0(\tF,K_{\tF})=5$.
 Since $h^0(D,K_D)=3$,
the image
$V$ of
$\HH{0}(\tF,K_{\tF})\to \HH{0}(D, {\cal O}_D(K_X+\tF))$ has dimension
$2$. On the other hand, $V$ contains the restriction to
$D$ of $\HH{0}(X,K_X)$,
that has dimension $3$ since $\phi$ maps $D$ to a conic. So
we have a contradiction, and we must conclude that if $X$ is of type
$I$, then $|F|$ ( and thus also $\phi^*|F|$) does not contain multiple
fibres.
\qed

\section{The examples} \label{esempi}
\setcounter{defn}{0}
\setcounter{equation}{0}
We describe here the known examples of surfaces $X$ such that the
canonical map of $X$ is $2$-$1$ on a canonically embedded surface
$\Si$, and we also present some new ones. We collect at the end of the
section some lemmas that are needed in the description of the examples.
\bigskip

\noindent {\bf $1.$  Examples with $X$ regular.}

The first example of a surface $X$ mapped non-birationally
onto a canonical surface by the canonical system was found
independently by several authors (\cite{bea}, \cite{babbage},
\cite{van}). One of the possible descriptions of the canonical image
$\Si$ is the following:
$\Si$ is a quintic surface in $\pp^3$, defined by the vanishing of
the determinant of a generic symmetric $5\times 5$ matrix $M$ of linear
forms. The singularities of $\Si$ are $20$ nodes, occurring precisely
where the rank of
$M$ drops by $2$, and they form an even set. The  double
cover of
$\Si$ branched over the nodes is a regular surface $X$ with $p_g=4$.
In \cite{ciro} p. 126 Ciliberto has remarked that the same
method can be used to produce similar examples, with
$\Si$ a canonical complete intersection in a projective space.
Notice that, if
$\Si$ is of type $(3,3)$, $(2,2,3)$ and
$(2,2,2,2)$, then the examples thus obtained are not on the Castelnuovo
line.
\bigskip

\noindent {\bf $2.$ Examples with $p_g(X)=5$ and $q(X)=2$.}

Surfaces $\Si$ with $p_g=5$ and $K^2=8$ are on the
Castelnuovo line and have been described in detail in \cite{ho}.  If
the canonical map is birational, then the canonical image is
isomorphic to the canonical model, and it is the intersection of a
quadric and a quartic in $\pp^4$. If the quadric is singular, then
$\Si$ is a Castelnuovo surface of type $(0,1,1)$ or $(0,0,2)$,
according to whether the quadric is singular at one point or along a
line. In the former case $\Si$ carries two different pencils of genus
$3$, while in the latter case $\Si$ is necessarily singular at the
points of intersection with the singular line of the quadric. (For a
generic choice of the quartic, these singularities will be $4$ nodes).

 Let now
$A$ be a principally polarized abelian surface,  let $K$ be the Kummer
surface of $A$ and let $q:A\to K$ be the projection onto the quotient.
If $D$ is a symmetric theta divisor, then the linear system $|2D|$ is
the pull-back from $K$ of a linear system
$|H|$. If $D$ is irreducible,  then $|H|$  embeds $K$ as a quartic
surface in
$\pp^3$; if $D$ is reducible, then $|H|$ maps $K$ $2$-$1$ onto the
smooth quadric in $\pp^3$. Let $f:\Si\to K$ be the double cover
branched on a curve $B$ of $|2H|$ not meeting the singular set of
$K$. If $B$ is smooth, then the singularities of
$\Si$ are $32$ nodes, which are the inverse image of the $16$
singular points of $K$. If $B$ has simple double points, then
$\Si$ has extra rational double points above the singularities of $B$.
 The sheaf
$f_*\OSi$ splits as
$\OO_K\oplus\OO_K(-H)$ and, using \ref{formuledoppi}, one computes:
$K^2_{\Si}=2H^2=8$,
$p_g(\Si)=p_g(K)+h^0(\OO_K(H))=5$,
$q(\Si)=h^1(\OO_K)+h^1(\OO_K(-H))=0$. Denote by
$\phi:X\to\Si$ the map obtained from $q:A\to K$ by base-changing
with $f$, and by
$\alpha:X\to A$ the map that completes the square as in diagram
\ref{diagram}. The map $\phi$ is branched precisely over the inverse
image of the singularities of $K$, while
$\alpha$ is  branched on $q^*B\in |4D|$, and $X$ is singular only above
the singularities of $q^*B$. One has:
$\alpha_*\OX=\OO_A\oplus\OO_A(-2D)$, and thus one may compute the
invariants of $X$ as above, and
obtain:
$p_g(X)=5=p_g(\Si)$, $K^2_X=16$, $q(X)=2$. (In fact, $\alpha$ is the
Albanese map of $X$). So the canonical map of
$X$ is the composition of $\phi$ with the canonical map of $\Si$.
If $D$ is irreducible, then
by our construction $\Si$ is isomorphic to the intersection in
$\pp^4$ of the cone over
$K$ with a quadric not passing through the vertex of the cone, and so
it is a canonical surface. As we have already explained at the
beginning, when the quadric is singular, namely when
$B$ is cut out on $K$ by a singular quadric of $\pp^3$, $\Si$ is a
Castelnuovo surface. In this case, it is easy to check that the genus
$3$ fibres are mapped to plane sections of $K$ by $\phi$, and that
their inverse images in $X$ are connected genus
$5$ curves; thus $X$ is a surface of type $I$.

Assume now that $D$ is reducible: then $A$ is isomorphic to the
product $E_1\times E_2$ of two elliptic curves, with
origins $O_i$ and, if
$\pi_i:A\to E_i$, $i=1,2$, are the projections, then
$D=\pi_1\inv(O_1)+\pi_2\inv(O_2)$. The map  $\pi_1\circ\alpha:X\to
E_1$, is an elliptic pencil of genus $3$ curves.  We wish to show
that, for a generic  choice of $B\in|2H|$, the generic fibres of
this pencil are not hyperelliptic. The subspace
$V=q^*\HH{0}(K,2H)\subset
\HH{0}(A,4D)$  is the subspace of even sections. It is possible to find
a basis
$\sigma^i_1,\ldots
\sigma^i_3,
\tau^i$ of $\HH{0}(E_i,4O_i)$, $i=1, 2$, such that the $\sigma^i_j$'s
are even and the $\tau^i$'s are
odd.   So
$V$ is spanned by the products
$\sigma^1_j\sigma^2_k$, $i,j=1,2,3$ and by $\tau^1\tau^2$. It follows
that the restriction of $V$ to a generic fibre of $\pi_1$  contains
sections that are not even. So,
for a generic choice of $B\in|2H|$, the inverse image in $X$ of a
generic fibre of
$\pi_1$ is not hyperelliptic by lemma \ref{hyperelliptic}.

 The maps
$\pi_1$ and $\pi_1\circ \alpha$ are compatible with the involutions on
$A$, on $X$ and on
$E_1$,
 and so they induce linear pencils $p_1:K\to
\pp^1$  and $p_1\circ
f:\Si\to\pp^1$. The generic fibre of
$p_1\circ f$ is the same as the generic fibre of $q_1\circ\alpha$, and
so it is  a non-hyperelliptic curve of genus
$3$. By lemma 1.1 of \cite{ak}, the canonical map of $\Si$ is not
composed with a pencil and it has degree
$\le 2$; on the other hand, the restriction of $|K_{\Si}|$ to a smooth
fibre  $F$ of $p_1\circ f$ is a subsystem of $|K_F|$. So we must
conclude that the restriction of the canonical map of $\Si$ to $F$ is
an embedding.
 Moreover, the  system $|K_{\Si}|$
contains the $f^*|2H|$, and so it separates the fibres of
$p_1\circ f$. So we conclude that the canonical map of $\Si$ is
birational and $\Si$ is a Castelnuovo surface (of type $(0,1,1)$). The
pull-back of the genus $3$ pencil
$p_1\circ f$ factors through the elliptic pencil $\pi_1\circ\alpha$,
thus it is not connected and $X$ is a surface of type $I\!I$.
\bigskip

\noindent {\bf $3.$ Surfaces of type $I\!I$ with $p_g=6$, $q=1$.}

From propositions \ref{constrII} and \ref{invariantsII}, it follows that
these examples arise from divisors $\Si$
of bidegree $(4,3)$ in $\pp^2\times \pp^1$ with only nodes as
singularities and having the following properties: 1)
the  pencil on $\Si$ induced by the projection
$p:\pp^2\times
\pp^1\to\pp^1$ has $4$ double fibres, 2) $\Si$ is smooth away from the
double fibres. Such  a surface $\Si$ has
$32$ nodes, $8$ on each of the double fibres, and these form an even
set. The double cover $\phi:X\to\Si$ branched over the nodes is a
surface of type $I\!I$.
In section \ref{conto} we produce explicitly such an example.
This enables us to describe a $16$-dimensional family
of non isomorphic surfaces $X$ of type $I\!I$ with the above invariants.
Let $U_1$ be the open subset of the $4$-fold product of $\pp^1$ with
itself consisting of the $4$-tuples $(z_1,z_2,z_3,z_4)$ such that
$z_i\ne z_j$ for $i\ne j$; let $U_2$ be the open subset of the
$4$-fold product of the space $\pp^5$ of conics with itself
consisting of the
$4$-tuples $(Q_1,Q_2,Q_3,Q_4)$ such that $Q_i$ is reduced,
$i=1,\ldots 4$, and
$Q_1^2$,
$Q_2^2$,
$Q_3^2$,
$Q_4^2$ represent independent points in the space
$\pp^{14}$ of quartics; finally, denote by $U_3\subset\pp^{55}$  the
open subset  of irreducible  divisors of bidegree $(4,3)$ in
$\pp^2\times\pp^1$, and let
$Z\subset U_1\times U_2\times U_3$ be the closed subset consisting of
the points
$(z_1,z_2,z_3,z_4;Q_1,Q_2,Q_3,Q_4;
\Si)$ such that the fibre of $\Si$ over $z_i$ is $Q_i^2$, for
$i=1,\ldots 4$. It is easy to check that the projection of $Z$ onto
$U_1\times U_2$ is surjective, and that the fibre of this projection
over a point of $U_1\times U_2$ is naturally isomorphic to $\pp^3$
minus the coordinate planes. So $Z$ is a smooth quasiprojective
variety of dimension $27$. Let now $U_0$ be the open
subset of
$Z$ consisting of the points such that the singularities of $\Si$ are
only nodes: the example of section \ref{conto} shows that $U_0$ is
nonempty. Moreover, the number $\nu(\Si)$ of nodes is a
lower-semicontinuous function of $\Si\in U_0$ and, by lemma
\ref{nodes},  one always has
$\nu(\Si)\ge 32$ . This minimum is attained in the example, and so
there is a nonempty open subset
$U\subset U_0$ such that $\Si$ has precisely $32$ nodes, occurring
on the double fibres. Notice that the restriction to $U$ of the
projection onto $U_3$ is a Galois cover  of its
image with Galois group $S_4$, the group action consisting simply in
changing the ordering of the double fibres of $\Si$. We  abuse
notation and also denote by
$U$ the image of $U$ in $U_3$. The double covers of surfaces $\Si\in
U$, branched over the nodes, form a $27$-dimensional family $W$ of
surfaces of type $I\!I$ with $q=1$ and $p_g=6$.
The group
$\pgl(2)\times \pgl(1)$ acts naturally on $U$, and thus on $W$. On the
other hand, it is easy to show that two surfaces
$X$,
$X'\in W$ are isomorphic iff they belong to the same orbit of the
action of
$\pgl(2)\times \pgl(1)$. Thus the geometric invariant theory quotient
of
$W$ by
$\pgl(2)\times
\pgl(1)$ is an irreducible variety of dimension $16$, parametrizing
non-isomorphic surfaces.
\bigskip

\noindent {\bf $4.$ An infinite family of surfaces of type $I$ with
q=2.} These examples are due to Beauville (see \cite{special}, $2.9$).
Let
$A$ be an abelian surface with an irreducible principal polarization
$D$, let
$K$ be the Kummer surface of $A$,
$q:A\to K$ the projection onto the quotient, and let
$|H|$ be the linear system on $K$ induced by $|2D|$. For
$n\ge 2$, consider a smooth hypersurface $G$ of bidegree $(1,n)$ in
$\pp^3\times
\pp^1$: the
projection onto $\pp^1$ exhibits $G$ as a $\pp^2$-bundle and the
hypersurfaces of bidegree $(1,k)$, $k\ge 0$ induce
effective tautological hyperplane sections of $G$.
Assume that $G$ intersects the singular locus of
$Y=K\times\pp^1$ transversely, and the intersection at smooth points
of $Y$ is transversal. (This certainly happens for a generic choice of
$G$). Then
the surface $\Si=G\cap Y$ has $16n$ nodes at the intersections with the
singular locus of $\Si$ and is smooth elsewhere.  Using adjunction on
$\pp^3\times\pp^1$, one sees that the hypersurfaces of bidegree
$(1, n-2)$ induce canonical curves of $\Si$;  so for $n\ge 3$ the
canonical system $|K_{\Si}|$ is very ample. It is not difficult to
check that the same is true for $n=2$. A straightforward computation
yields:
$p_g(\Si)=4n-3$,
$K^2_{\Si}=12n-16$. So $\Si$ is a Castelnuovo surface, and $G$, with
the natural projection onto
$\pp^1$, is isomorphic to the
$\pp^2$-bundle $\pabc$ containing it. Since the canonical divisor of
$\Si$ is induced by the hypersurfaces of bidegree $(1,n-2)$, one has
$a+b+c=p_g(\Si)-3=4n-6$, and $a\ge n-2$.

Now denote by
$X$ the pull-back of $\Si$ to $A\times\pp^1$ via the map $q\times 1$:
the surface
$X$ is smooth, the projection $X\to A$ is the Albanese map,  and
$q\times 1$ restricts to a double cover
$\phi:X\to\Si$ branched over
the nodes
$\Si$. Using adjunction on
$A\times\pp^1$ and
$K\times\pp^1$, one checks immediately that $p_g(X)=p_g(\Si)$. So
$\phi$ is the canonical map of $X$. Moreover, it is easy to see that
the inverse image of a fibre $F$ of the projection $\Si\to\pp^1$ is
connected, and thus $X$ is of type $I$.
\bigskip

The results that follow show that the surfaces $\Si$ of example
$4$ exist for all the admissible values of $a$, $b$ and
$c$.

\begin{lem}
Let $n\ge 3$ and $n-2\le a \le b\le c$ be integers such that
$a+b+c=4n-6$; then there exists a smooth divisor
$G\in\pp^3\times\pp^1$  of bidegree
$(1,n)$, and an isomorphism  $\pabc\to G$ such
that hypersurfaces of bidegree $(1,n-2)$ pull back to tautological
hyperplane sections of $\pabc$.

Write $4n-6=3\epsilon+\rho$, with $\rho$ and $\epsilon$ integers,
$0\le\rho <3$; a generic hypersurface $G$ of bidegree
$(1,n)$, with the polarization given by hypersurfaces of
bidegree $(1,n-2)$, is isomorphic to
$\pabc$ with the tautological hyperplane section, where
$a,b,c$ are as follows:

$\rho =0$, $a=b=c=\eps$;

$\rho =1$, $a=b=\eps$, $c=\eps+1$,

$\rho =1$, $a=\eps$, $b=c=\eps +1$.
\end{lem}
\proof
Let $T$ be the tautological hyperplane section and $L$ be the fibre
of the projection $p:\pabc\to\pp^1$; the divisor $T'=T-(n-2)L$ is base
point free, and the corresponding morphism $g:\pabc\to\pp^{n+2}$ is
birational. More precisely, if $a>n-2$ then $g$ is an embedding and if
$a=n-2$ then the image of $g$ is a cone over $\pp_{b,c}$. Let
$h:\pabc\to\pp^3$ be the morphism associated to a generic
$3$-dimensional subsystem of
$|T'|$:
$h$ has degree $n$ and maps the fibres of $\pabc$ linearly to planes
in
$\pp^3$. The morphism $h\times p:\pabc\to\pp^3\times\pp^1$ embeds
$\pabc$ as a divisor of type $(1,n)$, and hypersurfaces of bidegree
$(1, n-2)$ pull-back to elements of $|T|$ via $h\times p$.

To prove the second part of the statement, consider the
space $\pp^{4n+3}$ of divisors of bidegree $(1,n)$, and  the  dense
open subset
$U\subset
\pp^{4n+3}$  consisting of the smooth divisors.
If $k$ is an integer, then $h^0(G,\OO_G(1,k))$
is a lower semi-continuous function of $G\in U$. If, say,
$\rho=0$, then we have shown that there exists $G_0\in U$
such that $G_0$ with the polarization given by hypersurfaces of
bidegree $(1, n-2)$ is isomorphic to
$\pp_{\eps,\eps,\eps}$ with the tautological hyperplane sections.
This is equivalent to the condition
$h^0(G_0,\OO_{G_0}(1,n-\eps-3)=0$. Then, by semi-continuity, one
has
$h^0(G,\OO(1,n-\eps-3)=0$ on a dense open set $U_1\subset U$, and so
$G$ is isomorphic to $\pp_{\eps,\eps,\eps}$ for every $G\in U_1$. The
same argument shows the statement for $\rho=1,2$.
\qed

\begin{lem} Let $A$ be an abelian surface with an irreducible
principal polarization, let $K\subset \pp^3$ be the corresponding
Kummer quartic, and let $G$ be a smooth hypersurface of bidegree
$(1,n)$ in $\pp^3\times\pp^1$: there exists $\gamma\in \pgl(3)$
such that $\Si_{\gamma}= G\cap (\gamma K\times \pp^1)$ has $16n$ nodes,
occurring at the intersections of $G$ with the singular locus of
$K\times\pp^1$, and is smooth elsewhere.
\end{lem}
\proof
The proof consists simply in counting dimensions.
Let  $(\pp^3)^*$ be the space of planes in $\pp^3$, let
$K^*\subset(\pp^3)^*$ be the dual surface of $K$,  and let
$\psi:
\pp^1\to (\pp^3)^*$ be the map that associates to $z\in \pp^1$ the
plane $G\cap (\pp^3\times\{z\})$.
 We say that
$\gamma\in \pgl(3)$ is ``good'' if  $\gamma
K^*$ and
$\psi(\pp^1)$ intersect transversely at smooth points, and moreover
the intersection points are regular values of $\psi$ and do not lie
on the exceptional curves corresponding to the nodes of $K$. We are
going to show that if  $\gamma$ is ``good'', then it satisfies the
claim. Remark first of all that the points of the curve
$\psi(\pp^1)$ correspond to planes that are tangent to $\gamma K$ at
most at one smooth point. So  the surface
$\Si_{\gamma}$ has nodes at the points of intersections with the singular
set of
$\gamma K\times \pp^1$. To show that $\Si$ is smooth
elsewhere, notice that
$\Si_{\gamma}\cap(\pp^3\times\{z\})$ is just the intersection of
$\gamma K$ with the plane $\psi(z)$, and so
$\Si_{\gamma}$ can be singular only at points
$(x,z)\in G$ such that the plane
$\psi(z)$ is tangent to $\gamma K$ at $x$. A computation in local
coordinates shows that these points are also smooth if $\psi$ is
regular at $z$ and the curve
$\psi(\pp^1)$ meets $\gamma K^*$ transversely at $\psi(z)$.
In order to conclude the proof it is enough to remark that the
$\gamma$'s that are not ``good'' form a subset of dimension at most
$14$. This is a consequence of the following facts:
1) the subset of
$\pgl(3)$ consisting of the elements that map a  point
$x_1\in \pp^3$ to a point $x_2$ has dimension $12$,
2) the subset of
$\pgl(3)$ consisting of the elements that map a chosen point
$x_1\in \pp^3$  and a line $L_1$ through $x_1$ to a point
$x_2$ and a line $L_2$ through $x_2$ has dimension
$10$ .
\qed
The two previous lemmas together yield the following:
\begin{prop}
Let $A$ be an abelian surface with an irreducible
principal polarization, and let $a\le b\le c$ be integers such that
$a+b+c\equiv 2 (\mbox{mod}\,\,\, 4)$ and
$a\ge n:=(a+b+c+6)/4$; then there exist $X$ and $\Si$ as in example $4$
such that $\Si$ is a Castelnuovo surface of type $(a,b,c)$ and $A$ is
the Albanese variety of $X$.
\end{prop}
We close the section by proving the lemma needed in example $2$.
\begin{lem}\label{hyperelliptic}
Let $E$ be an elliptic curve with origin $O$; let $B\in |4O|$ be a
reduced divisor and let $f:C\to E$ be the double cover branched on
$B$, with ${\cal L}=2O$. Then $C$ is hyperelliptic if and only if $B$ is
symmetric with respect to the elliptic involution.
\end{lem}
\proof
 As it is explained in section
\ref{notation},
$C$ is isomorphic to a divisor $D\subset {\cal L}$; the line bundle
${\cal L}$ has a natural linearization and, if $B$ is symmetric, then  $D$
is easily seen to be also symmetric. Thus the involution on
${\cal L}$ induces an involution of $C$, whose fixed points are the inverse
images of points of order $2$ of $E$. Since $B$ is symmetric and
reduced, it does not contain any point of order $2$, and so the
involution has
$8$ fixed points on $C$. By the Hurwitz formula, the quotient of $C$ by
the involution is rational, and thus $C$ is hyperelliptic.

Conversely, assume that $C$ is hyperelliptic and denote by
$\phi:C\to\phi(C)$ the canonical map, with $\phi(C)$ a plane conic. If
$g:E\to\pp^1$ is the quotient map of the elliptic involution, then by
\ref{formuledoppi} the canonical system $\HH{0}(C,K_C)$ contains
$f^*\HH{0}(E, 20)= f^*g^*\HH{0}(\pp^1,\OO_{\pp^1}(1))$ as a subsystem.  So
one has a map
$\bar{f}:\phi(C)\to\pp^1$ such that the following diagram commutes:
\begin{equation}
\begin{array}{rcccl}
\phantom{1} & C &\stackrel{f}{\rightarrow} & E & \phantom{1} \\
\scriptstyle{\phi}\!\!\!\!\!\!\!\! & \downarrow & \phantom{1} &
\downarrow & \!\!\!\!\!\!
\scriptstyle{g}
\\
\phantom{1} & \phi(C) & \stackrel{\bar{f}}{\rightarrow} & \pp^1 &
\phantom{1}
\end{array}
\end{equation}
If we denote by $i_1$ the hyperelliptic involution on $C$ and by $i_2$
the elliptic involution on $E$, then it  follows immediately: $f\circ
i_1=i_2\circ f$.  In particular, if $R$ is the ramification divisor of
$f$ , then $i_1( R)=R$. Applying the Hurwitz formula to
$f$, one sees that
$R$  is a canonical divisor of
$C$; since $R$ is reduced, it contains no Weierstrass
point.  Thus, $R$ may be written as $x+i_1(x)+y+i_1(y)$, for some $x$,
$y\in C$ and $B=f(x)+i_2(f(x))+f(y)+i_2(f(y))$ is a symmetric divisor.
\qed

\section{Surfaces of type $I$ with $q(X)\ge 2$}\label{typeI}
\setcounter{defn}{0}
\setcounter{equation}{0}

In this section we study surfaces of type $I$ with $q\ge 2$ and we show
that they are all obtained as in example $4$ of section
\ref{esempi}. More precisely we prove the following:
\begin{thm}\label{main}
Let $\phi:X\to\Si$ be as in assumption \ref{ipotesi}; if $X$ is of
type $I$ and $q\ge 2$, then  $p_g(X)\equiv
1 (\mbox{mod}\,\,\, 3)$, $q(X)=2$, the Albanese surface $A$ of $X$ has
an irreducible principal polarization,  and $X$ can be
constructed as in example
$4$ of section \ref{esempi}, with $n=(p_g(X)+3)/4$.
\end{thm}
The proof of theorem \ref{main} requires some preliminary steps: first we
show that the Albanese variety $A$ of $X$  is a surface,  and then we prove
that $A$ is isomorphic to the Prym variety of the unramified double cover
$\tF\to F$, with $F$ a generic fibre. Thus $A$ is
principally polarized and  there is a map
$f:\Si\to K$, where $K$ is the Kummer quartic $K$ of $A$. Finally we show
that $f$ can be extended to a morphism $g:\pabc\to\pp^3$.
For the rest of the section, we will assume that
$\phi:X\to
\Si$ is as in assumption \ref{ipotesi}, that $X$ is of type $I$ and
that $q(X)\ge 2$. In particular, the pull-back to $X$ of the genus $3$
pencil
$|F|$ is a linear pencil $|\tF|$ of genus $5$.
\begin{lem}
The irregularity $q(X)$ is
equal to $2$.
\end{lem}
\proof
It suffices to show that $q(X)\le 2$.
 Notice that for a generic fibre
$F$, the restriction map $\HH{0}(\Si,K_{\Sigma}+F)\to
\HH{0}(F,K_F)$ is surjective, by the regularity of $\Sigma$. So
$\phi^*\HH{0}(\Si,K_{\Si}+F)\subset
\HH{0}(X,K_X+\tF)\to\HH{0}(\tF,K_{\tF})$
is a subspace whose image via the restriction map
$\HH{0}(X,K_X+\tF)\to\HH{0}(\tF,K_{\tF})$ has dimension $3$.
Since the pencil $|\tilde{F}|$ is linear, by Ramanujan vanishing
one has
$\HH{1}(X, K_X+\tF)=0$ and therefore the cokernel of the above
restriction map is isomorphic to
$\HH{1}(X,K_X)$. Since $\tF$ has genus $5$, it follows $q(X)\le 2$.
\qed

\begin{prop}\label{ppav} The Albanese variety $A$  of $X$ is a principally
polarized abelian surface, and the polarization $D$ of $A$ is irreducible.
\end{prop}

The above proposition is a consequence of the following lemmas, that
describe the Prym variety $Z$ of the cover $\tF\to F$, for a generic $F$
and show that $Z$ is naturally isomorphic to $A$.

 We start by reviewing
quickly the properties of Prym varieties that we need; for more
details and proofs the reader may consult chapter $12$ of \cite{lb}.
Let
$J$ be the Jacobian of $\tF$, and let $\gamma:\tF\to J$ be the period
map with base point $x_0\in \tF$. The Abel-Prym map with base point
$x_0$,
$\beta:\tF\to Z$, is defined as the composition $\hat{i}\circ \gamma$,
where $\hat{i}:J\to Z$ is a surjective morphism of abelian varieties
with connected kernel.
$Z$ is an abelian surface having a natural principal polarization $D$,
the restriction of
$\beta$ to
$\tF$ is an embedding and
the image of
$F$ is a divisor algebraically equivalent to $2D$, by Welters criterion
(\cite{lb}, page $373$).
\begin{lem}
The polarization $D$ on $Z$ is irreducible.
\end{lem}
\proof
Assume by contradiction that $D$ is reducible: then $Z$ is  a product
$E_1\times E_2$ of elliptic curves and $D$ is algebraically equivalent to
$\pi_1\inv(O_1)+\pi_2\inv(O_2)$, where $\pi_i$ is the projection onto
$E_i$, $i=1,2$, and  $O_i\in E_i$. For a suitable choice of $O_1$ and
$O_2$, the curve $\beta(\tF)$ in $E_1\times E_2$ is linearly equivalent to
$\pi_1\inv(2O_1)+\pi_2\inv(2O_2)$
 Denote by $j_i$ the
involution on $E_i$ that fixes $O_i$, $i=1,2$: then  $\beta(\tF)$is
invariant under
$j_1\times 1$ and
$1\times j_2$. The quotient of
$\beta(\tF)$ by the diagonal automorphism $j_1\times j_2$ is
isomorphic to $F$, and, via this isomorphism, $j_1\times 1$ induces an
involution of $F$ whose quotient is a plane section of the smooth
quadric in $\pp^3$. So $F$ is hyperelliptic, but this contradicts the
assumption that the canonical map of $\Si$ is birational.
\qed
Given any map $h:\tF\to Y$, with $Y$ a complex torus, there exist a
unique morphism of tori $\bar{\psi}:J\to Y$ and a unique translation
$\tau:Y\to Y$ such that
$\tau\circ h=\bar{\psi}\circ \gamma$. If, moreover, the map $h$ is
equivariant with respect to the $\Z/2$-actions given by the involution on
$\tF$ and by multiplication by $-1$ on
$Y$, then the kernel of
$\hat{i}$ is mapped to $0$ by $\bar{\psi}$; thus $\bar{\psi}$ induces
a morphism $\psi:Z\to Y$ such that $\tau\circ h=\bar{\psi}\circ
\beta$. We are interested in the case in which
$h$ is the restriction to $\tF$ of the Albanese map $\alpha:X\to A$.
Consistently with the above notation, we denote by $\bar{\psi}:J\to
A$ and $\psi:Z\to A$ the morphisms induced by the restriction of
$\alpha$.
\bigskip

\begin{lem}\label{iso}
The morphism $\psi:Z\to A$ is an isomorphism.
\end{lem}
By the above discussion, $\psi$ is an isomorphism iff $\bar{\psi}$ is
surjective and has connected kernel. In turn, if we consider the dual
map of $\bar{\psi}$, $\bar{\psi}^*:\Pic^0(X)\to\Pic^0(\tF)$, then the
above conditions are equivalent to $\bar{\psi}^*$ being injective. So
assume that
$\bar{\psi}^*$ is not injective and consider a torsion element
$\xi\in\ker\bar{\psi}^*$ of order $m>1$.  Let $r:X'\to X$ be the
unramified $\Z/m$-cover given by $\xi$: the restriction of $r$ to a
generic fibre $\tF$ is a disjoint union of $m$ components isomorphic
to $\tF$. Using the Stein factorization of the pull-back to $X'$ of
the pencil $\tF$, one gets the following commutative diagram, where
the vertical arrows are pencils with fibre $\tF$ and the horizontal
arrows are connected $\Z/m$-covers:
\begin{equation}
\begin{array}{rcccl}
\phantom{1} &X' &\stackrel{r}{\rightarrow} & X & \phantom{1} \\
\phantom{1} & \downarrow & \phantom{1} & \downarrow
&\phantom{1}
\\
\phantom{1} & B & \stackrel{\bar{r}}{\rightarrow} & \pp^1 & \phantom{1}
\end{array}
\end{equation}
The map $\bar{r}$ is ramified at at least $2$ points, while $r$,
which is obtained from $\bar{r}$ by base change and normalization, is
unramified. This implies that the fibres of the pencil $|\tF|$ over the
branch points of $\bar{r}$ are $m$-tuple fibres. But this contradicts
proposition \ref{nodes}.
\qed

 Finally, we put all the previous results together and get:
\smallskip

{\bf Proof of theorem \ref{main}:}
Consider the basic diagram \ref{diagram}: we wish to show that the map
$f:\Si\to K$ can be extended to a map $\bar{f}:\pabc\to\pp^3$ that
maps the fibres of $\pabc$ linearly to planes of $\pp^3$.
By lemma \ref{iso}, the Albanese variety $A$ can be identified with
the Prym variety $Z$ of a generic $\tF$. So, the fibres $\tF$ are
mapped to divisors of $|2D|$, where $D$ is the principal polarization
on $A$ (see proposition \ref{ppav}); as a consequence, the fibres
$F$ on $\Si$ are mapped isomorphically to plane sections of $K$ with
respect to the embedding as a quartic in $\pp^3$. If $F$ is a smooth
fibre of $\Si$, then there is a natural linear isomorphism between the
fibre of $\pabc$ containing $F$ and the plane in $\pp^3$ containing
$f(F)$. So we can define a rational map
$\bar{f}:\Si\to\pp^3$, such that its restriction to $\Si$ extends to
the morphism $f$. Let now $F_0$ be a fibre of $\pabc$ containing an
indeterminacy point of $\bar{f}$: the restriction of $\bar{f}$ to
$F_0$ is a degenerate projectivity, whose singular locus can either
be a point or a line. In the former case, the scheme theoretic
image via $f$ of the curve
$\Si\cap F_0$ would be a $4$--tuple line. This is impossible, because
a Kummer surface has no such plane section. If the indeterminacy
locus of $\bar{f}$ on $F_0$ were a line not contained in $\Si$, then
the curve $\Si\cap F_0$ would be contracted to a point, but this is
of course impossible. So the only possibility left is that the
indeterminacy locus of $\bar{f}$ on $F_0$ is a line $R$ contained in
$\Si$. Remark that every other component of $\Si\cap F_0$ is
contracted by $f$, and so the pull-back to $X$ of every component of
$\Si\cap F_0$ different from $R$ is contracted by $\alpha$. Arguing
as in the proof of proposition \ref{nodes}, one shows that $R$ can
contain at most $4$ nodes of $\Si$, so the pull-back $\tilde{R}$ of $R$
to
$X$ is a curve of genus $0$ or $1$. On the other hand,
the scheme-theoretic image $\Delta$  of the fibre of $|\tF|$
containing $\tilde{R}$ is supported on $\alpha(\tilde{R})$, but this is
impossible because $\Delta$ is an ample divisor on an abelian surface. So
we conclude that $\bar{f}$ is indeed a regular map. If we denote by
$p:\pabc\to\pp^1$ the projection map, then the map $\bar{f}\times
p:\pabc\to \pp^3\times
\pp^1$ embeds $\pabc$ as a divisor $G$ of bidegree $(1,n)$, for a
suitable value of $n$, and
$\Si$ is mapped isomorphically to $G\cap(K\times\pp^1)$. To determine
$n$, we use adjunction on $\pp^3\times\pp^1$ and
remark that divisors of bidegree $(1, n-2)$ cut out canonical curves
on
$\Si$, and so $p_g(X)=p_g(\Si)=a+b+c+3=4n-3\equiv 1 (\mbox{mod} 4)$,
$n=(p_g(X)+3)/4$. By lemma \ref{basic}, $\Si$ has $16n$ nodes,
occurring at the intersections of $G$ with the singular locus of
$K\times\pp^1$. So the intersection of $G$ with the singular locus of
$K\times\pp^1$ is transversal.
\qed

\section{Surfaces of type $I\!I$}\label{typeII}
\setcounter{defn}{0}
\setcounter{equation}{0}
In this section we describe surfaces of type $I\!I$ in detail and we
show that the invariants of these surfaces are bounded.

So here $X$ and $\Si$ are as in assumption \ref{ipotesi}, and
moreover the pull-back $\tF$ of a generic $F$ is
disconnected. The Stein factorization of the pencil $|\tF|$ gives rise
to the following commutative diagram, where $p$ denotes
the pencil $|F|$ and $\tilde{p}$ denotes the connected fibration on $X$
through which
$|\tF|$ factors:
\begin{equation}\label{diagramII}
\begin{array}{rcccl}
\phantom{1} & X & \stackrel{\phi}{\rightarrow} & \Si \\
\scriptstyle{\tilde{p}}\!\!\!\!&\downarrow &\phantom{1} & \downarrow &
\!\!\!\!\scriptstyle{p}
\\
\phantom{1} & B & \stackrel{\bar{\phi}}{\rightarrow} & \pp^1
&\phantom{1}
\end{array}
\end{equation}
The curve $B$ is smooth and the map $\bar{\phi}$ is a double cover. We
introduce a new invariant of
$X$, the genus $g$ of $B$. Notice that $g\le q(X)$
\begin{thm}\label{constrII}
The surface $\Si$ has precisely $2g+2$ double fibres, occurring at
the branch points of $\bar{\phi}$, and $c\le g$. Conversely, if $\Si$ a
Castelnuovo surface of type $(a, b,c)$ with only nodes as
singularities, with $c\le g$, having $2g+2$ double fibres, and smooth
outside the double fibres, then $\Si$ has $16g+16$ nodes which  form
an even set, and the double cover $\phi:X\to\Si$ branched over the
nodes is a surface of type $I\!I$.
\end{thm}
\proof
For  $i=1,\ldots 2g+2$, denote by $x_i$ the ramification points
of $\bar{\phi}$, by
$y_i\in \pp^1$ the image of $x_i$, and by $F_i$ and $\tF_i$ the fibres
of
$p$ and $\tilde{p}$ over $y_i$ and $x_i$ respectively. In diagram
\ref{diagramII}, the map
$\phi$ is obtained from
$\bar{\phi}$ by base and normalization: since $\phi$ is
unramified in codimension $1$, the
$F_i$'s are double fibres, $i=1,\ldots 2g+2$,and they contain all the
nodes of $\Si$. So, by proposition \ref{nodes}, $\Si$ has precisely
$16g+16$ nodes. In order to show that $c\ge g$, we construct
explicitly $X$ as the normalization of a divisor in a $\pp^2$-bundle.
Denote by $\tpabc$ the pull-back of $\pabc$ to $B$; so
$\tpabc=\Proj(\bar{\phi}^*(\OP(a)\oplus\bar{\phi}^*(\OP(b)
\oplus\bar{\phi}^*(\OP(c)$. If $T$ and $L$ are the tautological
hyperplane section and a fibre of  $\pabc$, and $\tilde{T}$, $\tilde{L}$,
are the pull-backs of $T$ and $L$ to
$\tpabc$,  then the fibre product
$W$ of
$p$ and
$\bar{\phi}$ is a divisor in $\tpabc$ linearly equivalent to
$4\tilde{T}\!-\!(a+b+c-2)\tilde{L}$. The singular locus of $W$ consists of
$2g+2$ double curves, that are  the intersections of $W$ with
the fibres of $\tpabc$ over $x_1,\ldots x_{2g+2}$.  One has:
$K_{\tpabc}=-3\tilde{T}+\tilde{p}^*(K_B)+(a+b+c)\tilde{L}$ and
$K_{\tpabc}+W=\tilde{T}+\tilde{p}^*(K_B)$. So the
canonical curves of
$X$ correspond to sections of
$\tilde{T}+\tilde{p}^*(K_B+\bar{\phi}^*(\OP(2))$ vanishing on the double
curves of $W$, namely to sections of $\tilde{T}+tilde{p}^*(K_B-(x_1+\ldots
+x_{2g+2}+\bar{\phi}^*(\OP(2))$. The Hurwitz formula shows that $K_B$ is
linearly equivalent to
$x_1+\ldots+ x_{2g+2}-\bar{\phi}^*(\OP(2)$, and so the canonical system of
$X$ is the pull-back of
$\HH{0}(\tpabc,\tilde{T})$.
By assumption, we have
$p_g(X)=p_g(\Si)=a+b+c+3$; on the other hand,
$h^0(\tpabc,\tilde{T})=h^0(B,\bar{\phi}^*(\OP(a))+
h^0(B,\bar{\phi}^*(\OP(b)) +h^0(B,\bar{\phi}^*(\OP(c))$. Applying
\ref{formuledoppi} to the double cover
$\bar{\phi}$ (with ${\cal L}=\OP(g+1)$) yields for any integer
$k\ge 0$:
$h^0(B,\bar{\phi}^*\OP(k))=h^0(\pp^1,\OP(k))+h^0(\pp^1,\OP(k-g-1))=
k+1+h^0(\pp^1,\OP(k-g-1))$. So it follows that $c\le g$.

Conversely, assume that $\Si$ is a divisor in $\pabc$ linearly
equivalent to $4T-(a+b+c-2)L$, with only nodes as singularities, with
exactly $2g+2$ double fibres, $c\le g$, and assume that $\Si$ is smooth
away from the double fibres. Then by proposition
\ref{nodes} $\Si$ has $16g+16$ nodes. Let $D\subset \pabc$ be the
sum of the fibres of $p:\pabc\to\pp^1$ that are
double for $\Si\to\pp^1$: $D$ is smooth,
it is linearly equivalent to $2D'$, where $D'=(g+1)F$, it contains all
the nodes of $\Si$ and, finally, the restriction of $D$ to $\Si$ is a
union of double curves. Therefore, by proposition \ref{criterio}, the
nodes of $\Si$ form an even set. If $\phi:X\to\Si$ is the double cover
branched over the nodes, then the above computations show that $X$
is a smooth surface such that $p_g(X)=p_g(\Si)$. Let
$\epsilon:S\to\Si$ be the minimal resolution of the singularities of
$\Si$, let $E$ be the exceptional divisor of $\epsilon$, let $Z$ be the
sum of the strict transforms of supports of the double fibres of $\Si$ and
let
$\tilde{\phi}:\tilde{X}\to S$ be obtained from
$\phi$ by base change with $\epsilon$: $\tilde{\phi}$ is a smooth
double cover branched on $E$, and the
line bundle
${\cal L}$ associated with the cover is  equal $(g+1)\epsilon^*F-Z$. 	The
restriction of ${\cal L}$ to a generic fibre $F$ is trivial,  so the
inverse image in
$\tilde{X}$ (and in
$X$) of a generic $F$ is disconnected, and $X$ is of type $I\! I$.
\qed

We will now give bounds for the invariants of $X$
and relations between them.
\begin{lem}\label{q-g}
The following relations between the invariants of $X$ hold:
$$0\le q-g=3g+3-p_g.$$
\end{lem}
\proof The inequality $q-g\ge 0$ is a consequence of the fact that
there exists a dominant map $\tilde{p}:X\to B$, with $B$ a curve of
genus
$g$. The equality
$q-g=3g+3-p_g$ is equivalent to the formula \ref{nodes}, since by
propositions \ref{constrII} $\Si$ has
 $16g+16$ nodes.
\qed
\begin{prop}\label{invariantsII}
The numerical possibilities for the invariants of $X$ are the
following:
\smallskip

\noindent a) $p_g=3g+3$,\quad $q=g$,\quad\quad $a=b=c=g$,\quad $0<g\le
26$;

\noindent b) $p_g=3g+2$,\quad $q=g+1$, \quad $a=g-1$, $b=c=g$, \quad
$0<g\le 17$;

\noindent c) $p_g=3g+1$,\quad $q=g+2$,\,\,\,\quad $a=b=g-1$, $c=g$ or
 $a=g-2$,
$b=c=g$, \quad
$0<g\le 8$.
\end{prop}
\proof
The topological Euler characteristic of the minimal desingularization
$S$ off $\Si$ can be computed from Noether's formula as follows:
$$c_2(S)=12\chi(\OO_S)-K^2_S=9p_g(S)+19.$$
On the other hand the following formula (see \cite{bpv}, page 97),
in which
$e(D)$ represents the topological characteristic of the support of a
divisor
$D$, expresses
$c_2(S)$ in terms of the base and of the singular fibres of the
pull-back to $S$ of the fibration
$F$ on $\Si$, that we also denote by $F$:
$$c_2(S)=e(\pp^1)e(F)+\sum_{F'\,
singular}e(F')-e(F).$$
(The term
$e(F')-e(F)$ is always non-negative, see \cite{bpv}, page 97). From
proposition
\ref{nodes}, it follows that if $F'$ is the pull-back of a double
conic on $\Si$, then
$e(F)=10$ or $11$, according to whether the conic is smooth or not. So,
recalling that by proposition
\ref{constrII} there are $2g+2$ double fibres on $\Si$, and comparing
the two expressions for
$c_2(S)$  one obtains: $c_2(S)\ge 2(-4) +(2g+2)14$, namely
$9p_g\ge 28g+1$. The statement now follows in view of lemma \ref{q-g}.
\qed
\begin{rem}
Notice that at least for $g=1$ possibilities a) and b) actually
occur, as it is shown by examples $3$ and $2$ of section
\ref{esempi}. We do not know examples for possibility c).
\end{rem}

\section{Appendix: a computation with Macaulay}\label{conto}

We describe here how we have used Macaulay  (\cite{macaulay}) to show the
existence of a divisor $\Si$ of bidegree $(4,3)$ in $\pp^2\times\pp^1$
with the properties required in example $3$ of section \ref{esempi}.
We use the notation introduced there.

We consider homogeneous coordinates $(s,t)$ in $\pp^1$ and $(x_0,x_1,x_2)$
in
$\pp^2$ and set: $z_1=(1,0)$, $z_2=(0,1)$, $z_3=(1,-1)$,
$z_4=(1,1)$;\quad
$Q_1=x_0^2+x_1^2+x_2^2$, $Q_2=x_0^2+x_1^2-x_2^2$,
$Q_3=x_0x_1+x_0x_2+x_1x_2$,
$Q_4=x_0^2+5x_0x_1+7x_0x_2+2x_1^2+11x_1x_2+3x_2^2$. We start by writing
down the equation $h$ of $\Si$:

\begin{verbatim}
Macaulay version 3.0, created 12 September 1994

1% ring R ! characteristic (if not 31991)       ?
! number of variables  ? 5
! 5 variables, please ? stx[0]-x[2]
! variable weights (if not all 1)   ?
! monomial order (if not rev. lex.)   ?
largest degree of a monomial        : 217
1% ideal f1
! number of generators ? 1
! (1,1) ? x[0]2+x[1]2+x[2]2
1% ideal f2
! number of generators ? 1
! (1,1) ? x[0]2+x[1]2-x[2]2
1% ideal f3 ! number of generators ? 1
! (1,1) ? x[0]x[1]+x[0]x[2]+x[1]x[2]
1% ideal f4 ! number of generators ? 1
! (1,1) ? x[0]2+5x[0]x[1]+2x[1]2+7x[0]x[2]+11x[1]x[2]+3x[2]2
1% poly s1 (s-t)*(s+t)*s
1% poly s2 (s-t)*(s+t)*t
1% poly s3 (s-t)*st
1% poly s4 (s+t)*st
1% mult f1 f1 g1
1% mult f2 f2 g2
1% mult f3 f3 g3
1% mult f4 f4 g4
1% mult s1 g1 h1
1% mult s2 g2 h2
1% mult s3 g3 h3
1% mult s4 g4 h4
1% add h1 h2 h
1% add h h3 h
1% add h h4 h
%
\end{verbatim}
Next we show that the singularities of $\Si$ are at most nodes. This is
a local computation,  that has to be repeated fo each of the $6$ standard
open affine subsets.
Consider for instance $U=\{sx_0\ne 0\}\subset\pp^1\times\pp^2$: we
identify $U$ with $\A^3\subset\pp^3$ and then consider the closure in
$\pp^3$ of
$\Si\cap U$, defined by the equation $h_{s0}$. The plane at infinity is
 $w=0$. The ideal
$I$ of the locus in
$\pp^3$ of the singular points of $h_{s0}=0$ that are not nodes is
generated by the derivatives of $h_{s0}$ and by the
$3\times 3$ minors of the Hessian matrix  $hh_{s0}$ of $h_{s0}$.
Computing the standard basis of $I$ and using it to  reduce $w^{30}$ one
gets $0$, namely $w^{30}\in I$ and therefore the singularities of
$\Si\cap U$ are at most nodes. Here is the transcript of the Macaulay
session (slightly edited):
\begin{verbatim}
1% ring Ss0
! characteristic (if not 31991)  ?
! number of variables                 ? 4
!   4 variables, please               ? tx[1]x[2]w
! variable weights (if not all 1)     ?
! monomial order (if not rev. lex.)   ?
  largest degree of a monomial        : 512
1% rmap fs0 R Ss0
! s ---> ? w
! t ---> ? t
! x[0] ---> ? w
! x[1] ---> ?x[1]
! x[2] ---> ? x[2]
1% ev fs0 h hs0
1% setring Ss0
1% jacob hs0 jhs0
1% jacob jhs0 hhs0
1% flatten jhs0 jhs0
1% wedge hhs0 3 whs0
[189k][252k]
1% flatten whs0 whs0 [315k]
1% concat whs0 jhs0
1% lift-std whs0 whs0std
6.7.8.[378k]9.10.11.[441k]12.[504k][567k]13.14.15.[630k][692k]16.
[755k][818k][881k]17.[944k][1007k]18.[1070k]19.20.[1133k]21.
[1196k] computation complete after degree 21
% ideal w
! number of generators ? 1
! (1,1) ? w30
% reduce whs0std w red
% type red 0
\end{verbatim}
The computation in the other $5$ affine open sets goes exactly in the
same way. Now, to finish the computation it is enough to show that the
singular locus of $\Si$, that we already know to be reduced and of
dimension $0$, has length $32$. In fact, by prop. \ref{nodes}, each of
the $4$ double fibres contains $8$ nodes and therefore $\Si$ is smooth
outside the double fibres.
First one embeds $\Si$ in $\pp^5$ via the Segre embedding:
\begin{verbatim}
1% ring S ! characteristic (if not 31991)       ?
! number of variables                 ? 11 [126k]
!  11 variables, please               ? stx[0]-x[2]y[0]-y[5]
! variable weights (if not all 1)     ? 1:5 2:6
! monomial order (if not rev. lex.)   ? 5 1 1 1 1 1 1
largest degree of a monomial      : 217 512 512 512 512 512 512
1% fetch h h
1% ideal j
! number of generators ? 6
! (1,1) ? sx[0]-y[0]
! (1,2) ? sx[1]-y[1]
! (1,3) ? sx[2]-y[2]
! (1,4) ? tx[0]-y[3]
! (1,5) ? tx[1]-y[4]
! (1,6) ? tx[2]-y[5]
1% concat h j
1% std h hst
23.4.5.6.7.8.9.10.11.12.13.14.[189k]15.16.17.
computation complete after degree 17
1% elim hst helim
1% ring R
! characteristic (if not 31991)       ?
! number of variables                 ? 6
!   6 variables, please               ? y[0]-y[5]
! variable weights (if not all 1)     ?
! monomial order (if not rev. lex.)   ?
  largest degree of a monomial        : 117
1% fetch helim h
\end{verbatim}
Now $h$ is the ideal of  $\Si$ in $\pp^5$; the singular locus of $\Si$
is defined by the equations of $\Si$ and by the $3\times 3$ minors of
the Jacobian matrix of the equations of $\Si$:
\begin{verbatim}
1% std h hst 23.4.5.6.7.8. computation complete after degree 8
1% jacob hst jh
1% wedge jh 3 sing
[252k][315k][378k][441k][504k][567k][630k]
1% flatten sing sing
[692k][755k][818k][881k][944k]
1% concat sing hst
1% std sing singst 0123.4.5.6.7.8. computation complete after
 degree 8
1% degree singst
codimension : 5
degree      : 32
\end{verbatim}
The last line shows that the singularities of $\Si$ are $32$, as required.

\end{document}